\documentclass[11pt]{amsart}

\usepackage{amssymb,amsmath}
\usepackage{amsthm}
\usepackage{qsymbols}

\newtheorem{thm}{Theorem}
\newtheorem{lem}[thm]{Lemma}
\newtheorem{prop}[thm]{Proposition}
\newtheorem{cor}[thm]{Corollary}
\newtheorem{defi}[thm]{Definition}
\newtheorem{rem}[thm]{Remark}

\newcommand{\X}{\mathcal{X}}
\newcommand{\Y}{\mathcal{Y}}
\newcommand{\PX}{\mathcal{P}(\X)}

\newcommand{\dom}{d_\omega(\,.\,,\,.\,)}
\newcommand{\domxy}{d_\omega(x,y)}

\newcommand{\la}{\left\{}
\newcommand{\ra}{\right\}}
\newcommand{\lp}{\left(}
\newcommand{\rp}{\right)}
\newcommand{\lc}{\left[}
\newcommand{\rc}{\right]}

\newcommand{\R}{\mathbb{R}}
\newcommand{\N}{\mathbb{N}}
\newcommand{\Rd}{\mathbb{R}^d}
\newcommand{\Rdn}{\lp\mathbb{R}^d\rp^n }

\def\AArm{\fam0 \rm}%
\newdimen\AAdi%
\newbox\AAbo%
\def\AAk#1#2{\setbox\AAbo=\hbox{#2}\AAdi=\wd\AAbo\kern#1\AAdi{}}%
\newcommand{\1}{{\ensuremath{{\AArm 1\AAk{-.8}{I}I}}}}

\newcommand{\E}{\mathbb{E}}

\newcommand{\Hnm}{\operatorname{H}(\nu\mid\mu)}
\newcommand{\ent}{\operatorname{Ent}}
\newcommand{\Var}{\operatorname{Var}}

\newcommand{\Capa}{\operatorname{Cap}}

\newcommand{\Id}{\operatorname{Id}}
\newcommand{\sgn}{\operatorname{sgn}}

\newcommand{\Tnm}{\mathcal{T}_c(\nu,\mu)}

\renewcommand{\a}{\alpha}

\newcommand{\T}{\mathbb{T}}
\newcommand{\SG}{\mathbb{SG}}
\newcommand{\LO}{\mathbb{LO}}
\newcommand{\BLO}{\mathbb{BLO}}

\begin{document}
\title[Poincar\'e inequalities for non euclidean metrics \ldots]{Poincar\'e inequalities for non euclidean metrics and transportation cost inequalities on $\R^d.$}
\author{Nathael Gozlan}

\date{\today}

\address{Universit\'e Paris-Est - Laboratoire d'Analyse et de Math\'ematiques Appliqu\'ees (UMR CNRS 8050), 5 bd Descartes, 77454 Marne la Vall\'ee Cedex 2, France}
\email{nathael.gozlan@univ-mlv.fr}

\keywords{Poincar\'e inequality, Transportation cost inequalities, Concentration of measure, Logarithmic-Sobolev inequalities}
 \subjclass{60E15 \and 26D10}
\maketitle
\begin{center}
 \textsc{Universit\'e Paris-Est,}\\
 Laboratoire Analyse et Math\'ematiques Appliqu\'ees,\\
 UMR CNRS 8050,\\
 5 bd Descartes, 77454 Marne la Vall\'ee Cedex 2, France
 \end{center}
\begin{abstract}
In this paper, we consider Poincar\'e inequalities for non euclidean metrics on $\R^d$. These inequalities enable us to derive precise dimension free concentration inequalities for product measures. This technique is appropriate for a large scope of
concentration rate: between exponential and gaussian and beyond.
We give different equivalent functional forms of these Poincar\'e type inequalities in terms of transportation-cost inequalities and infimum convolution inequalities. Workable sufficient conditions are given and a comparison is made with generalized Beckner-Latala-Oleszkiewicz inequalities.
\end{abstract}

\section{Introduction}
\subsection{Poincar\'e inequality and concentration of measure}
One says that a probability measure on a metric space $(\X,d)$ satisfies a Poincar\'e inequality also called spectral gap inequality with the constant $C$, if for all locally Lipschitz function $f$, one has
\begin{equation}\label{Poincaré}
\Var_{\mu}(f)\leq C\int |\nabla f|^2\,d\mu,
\end{equation}
where the length of the gradient is defined by
\begin{equation}\label{length}
|\nabla f|(x):=\limsup_{y\to x}\frac{|f(x)-f(y)|}{d(x,y)}.
\end{equation}
(when $x$ is not an accumulation point of $\X$, one defines $|\nabla f|(x)=0$).

It is well known since the works \cite{Gromov1983}, \cite{Aida1994},\cite{Aida1994a} and \cite{BL97} that the inequality (\ref{Poincaré}) implies dimension free concentration inequalities for the product measures $\mu^n$, $n\geq 1$. For example, in \cite{BL97}, M. Ledoux and S.G. Bobkov proved the following theorem (see \cite[Corollary 3.2]{BL97})

\begin{thm}[Bobkov-Ledoux]\label{BL1}
If $\mu$ satisfies (\ref{Poincaré}), then for every bounded function $f$ on $\X^n$ such that $\displaystyle{\sum_{i=1}^n |\nabla_i f|^2\leq a^2}$ and $\displaystyle{\max_{i=1,\ldots,n}|\nabla_i f|\leq b}$, $\mu^n$ a.e. (where $|\nabla_if|$ denotes the length of the gradient with respect to ith the coordinate) one has
\begin{equation}\label{deviation}\forall t\geq 0,\qquad \mu^n\lp f\geq \int f\,d\mu^n + t\rp\leq \exp\lp-\min\lp \frac{t^2}{C\kappa^2a^2},\frac{t}{\sqrt{C}\kappa b}\rp\rp,\qquad \text{with } \kappa=\sqrt{18e^{\sqrt{5}}}.
\end{equation}
\end{thm}

Another way to express the concentration of the product measure $\mu^n$ is the following:
\begin{cor}[Bobkov-Ledoux]\label{Cor-concentration}
Let $\mu$ be a probability measure on $\X$ satisfying the Poincar\'e inequality (\ref{Poincaré}) on $(\X,d)$ with the constant $C>0$. Define $K(C)=\alpha(\frac{1}{\sqrt{C}\kappa})/16$, where as before $\kappa=\sqrt{18e^{\sqrt{5}}}$, then for all subset $A$ of $\X^n$ with $\mu^n(A)\geq 1/2$,
\begin{equation}\label{concentration}
\forall h\geq 0,\qquad \mu^n\lp A^h\rp\geq 1-e^{-K(C)h},
\end{equation}
where the set $A^h$ is the enlargement of $A$ defined by \[A^h=\la y\in\X^n  : \inf_{x\in A}\sum_{i=1}^n \a(d(x_i,y_i))\leq h\ra,\]
where $\a(u)=\min(|u|,u^2)$ for all $u\in \R$.
\end{cor}
The inequality (\ref{concentration}) can be easily derived from Theorem \ref{BL1} (see \cite{BL97} or Section \ref{sec concentration} of the present paper). Inequalities such as (\ref{concentration}) were first obtained by M. Talagrand in different articles using completely different techniques (see e.g. \cite{Talagrand1995}).

If $\mu$ satisfies (\ref{Poincaré}) on $\Rd$ equipped with its standard euclidean norm $|\,\cdot\,|_2$, then (\ref{concentration}) can be rewritten in a more pleasant way: for all subset $A$ of $\Rdn$ with $\mu^n(A)\geq 1/2$,
\begin{equation}\label{concentration2}
\forall h\geq0,\qquad \mu^n\lp A +\sqrt{h}B_2 + hB_1\rp\geq 1-e^{-hK(C)}
\end{equation}
with the same constant $K(C)$ as above.
The archetypic example of a measure satisfying (\ref{Poincaré}) is the exponential measure on $\R^d$ $\nu_1^d$, where $d\nu_1(x)=\frac{1}{2}e^{-|x|}\,dx$. For this probability, (\ref{concentration2}) cannot be improved (a version of (\ref{concentration2}) with sharp constants has been established by Talagrand in \cite{Talagrand1991} see also Maurey \cite[Corollary 1]{Mau91}). Thus (\ref{concentration2}) expresses that the probability measures $\mu^n$ concentrate at least as fast as the exponential measure on $\Rdn$.

Some probability measures concentrate faster than the exponential measure. For example, the standard gaussian measure $\gamma^m$ on $\R^m$ verifies for all $A\subset \R^m$ with $\gamma^m(A)\geq 1/2$,
\begin{equation}\label{gaussian}
\forall h\geq 0,\qquad \gamma^m(A+hB_2)\geq 1-e^{-h^2/2}.
\end{equation}
There is absolutely no hope to derive such a bound from the classical Poincar\'e inequality (\ref{Poincaré}) on $\R^m$ equipped with the euclidean norm. The inequality (\ref{gaussian}) requires other tools. For example (\ref{gaussian}) follows from the  Logarithmic Sobolev inequality, introduced by L. Gross in \cite{Gross1975}, which is strictly stronger than (\ref{Poincaré}) (see \cite[Chapter 5]{Led}).

\subsection{Changing the metric improves the concentration}The aim of this paper is to show that replacing in (\ref{Poincaré}) and (\ref{length}) the standard euclidean norm $|\,.\,|_2$ by another metric makes possible to reach a large scope of concentration properties including gaussian or even stronger behaviors. The metrics we are going to equip $\R^d$ with are of the form:
\begin{equation}\label{metric}
\forall x,y\in \R^d,\qquad d_\omega(x,y)=\left[\sum_{i=1}^d |\omega(x_i)-\omega(y_i)|^2\right]^{1/2},
\end{equation}
where, in all the paper, we will assume that $\omega:\R\to\R$ verifies:
\begin{itemize}
\item[$\bullet$] $\omega$ is such that $x\mapsto \omega(x)/x$ is non decreasing on $(0,+\infty),$
\item[$\bullet$] $\omega$ is non negative on $\R^+$,
\item[$\bullet$] $\omega$ is such that $\omega(-x)=-\omega(x)$, for all $x\in \R$.
\end{itemize}
Note that the first assumption is verified as soon as $\omega$ is convex on $\R^+$ with $\omega(0)=0$.
\begin{defi}
One says that a probability measure $\mu$ on $\R^d$ satisfies the inequality $\SG(\omega,C)$ (resp. $\SG(C)$) if $\mu$ satisfies the Poincar\'e inequality (\ref{Poincaré}) for the distance $\dom$ defined by (\ref{metric}) (resp. for the standard euclidean metric) with the constant $C>0$.
\end{defi}

Let us give a first example:
\begin{prop}\label{concentration omegap}
Let $\omega_p(x)=\max(x,x^p)$ on $\R^+$ with $\omega_p(-x)=-\omega_p(x)$ for all $x \in \R$.\\
Suppose that $\mu$ satisfies the inequality $\SG(\omega_p,C)$ on $\R^d$ for some $C>0$.\\

If $p\in [1,2]$, then for all $n\geq 1$ and all $A\subset \Rdn$,
\[\forall h\geq 0,\quad \mu^n\lp A + 2\sqrt{h} B_2 + 2h^{1/p}B_p\rp\geq 1 - e^{-K(C)h/d}.\]

If $p\geq 2$, then for all $n\geq 1$ and all $A\subset \Rdn$,
\begin{align*}
\forall h\geq 0,&\quad \mu^n\lp A + 2\sqrt{h} B_2\rp\geq 1 - e^{-K(C)h/d},\\
\text{and}\quad\forall h\geq 0,&\quad \mu^n\lp A + 2h^{1/p}B_p\rp\geq 1 - e^{-K(C)h/d}.
\end{align*}
(where $K(C)$ is defined in Corollary \ref{Cor-concentration})
\end{prop}

This result will be easily deduced from Corollary \ref{Cor-concentration} and from an elementary comparison between the metric $d_{\omega_p}(\,.\,,\,.\,)$ and the norms $|\,.\,|_p$.
In particular, it will follow from our general sufficient conditions that, for all $p\in [1,+\infty)$, the probability measure $d\nu_p(x)=\frac{1}{Z_p}e^{-|x|^p}\,dx$ verifies $\SG(\omega_p,C)$ for some $C$. The interest of our approach is to give a somewhat unified picture of the concentration of measure phenomenon.

\subsection{Presentation of the results}Before going into further details in the presentation of our results, let us introduce some notations and conventions.

\subsubsection{Notations} The map $\omega$ is defined on $\R$ but we will also denote by $\omega$ the map defined on $\R^m$ (for every $m\geq 1$) by $(x_1,\ldots,x_m)\mapsto (\omega(x_1),\ldots,\omega(x_n))$.
The image of a probability measure $\mu$ on a space $\X$ under a measurable map $T:\X\to\Y$ will be denoted by $T^\sharp\mu$. We recall that is is defined by
\[\forall A\subset \Y,\qquad T^\sharp\mu(B)=\mu \lp T^{-1}(A)\rp.\]

\subsubsection{Links with the classical Poincar\'e inequality} In Section \ref{Sec SG}, we prove the concentration results and we relate the exotic Poincar\'e inequalities $\SG(\omega,C)$ to (weighted) forms of the classical Poincar\'e inequality:
\begin{prop}\label{equivalence}
Let $\mu$ be a probability measure on $\Rd$ and $C$ a positive number. The following properties are equivalent.
\begin{itemize}
\item[(i)] The probability measure $\mu$ verifies $\SG(\omega,C).$
\item[(ii)] The probability measure $\omega^\sharp\mu$ verifies $\SG(C).$
\item[(iii)] The probability measure $\mu$ satisfies the following weighted Poincar\'e inequality:
\begin{equation}\label{weightedPoincaré}
\forall f,\qquad \Var_\mu(f)\leq C\int \sum_{i=1}^d \frac{1}{\omega'(x_i)^2}\lp\frac{\partial f}{\partial x_i}(x)\rp^2\,d\mu(x),
\end{equation}
for all $f:\Rd\to\R$ such that $f\circ \omega^{-1}$ is of class $C^1.$
\end{itemize}
\end{prop}
Observe that this proposition furnishes a huge collection of examples. Indeed, with a slight abuse of notations, one has
\[\lp\omega^{-1}\rp^\sharp \SG(C)\subset \SG(\omega, C).\]

\subsubsection{Sufficient conditions for $\SG(\omega,C)$} In Section \ref{Sec CS} we addressed the problem of finding workable sufficient conditions for Poincar\'e inequalities $\SG(\omega,C)$. The strategy is dictated by Proposition \ref{equivalence}. Namely, a probability $\mu$ satisfies $\SG(\omega,C)$, if and only if the measure $\omega^\sharp\mu$ satisfies $\SG(C)$. So all we have to do is to apply to the measure $\omega^\sharp \mu$ one of the known criteria for the classical Poincar\'e inequality.

In dimension one, one has a necessary and sufficient condition for $\SG(\omega,C)$:
\begin{prop}\label{CSdim1}
A probability measure $\mu$ on $\R$ absolutely continuous with density $h>0$ satisfies the inequality $\SG(\omega,C)$ for some $C>0$ if and only if
\begin{equation}
D_\omega^-=\sup_{x\leq m}\mu(-\infty,x]\int_x^m\frac{\omega'(u)^2}{h(u)}\,du<+\infty \quad\text{and}\quad D_\omega^+=\sup_{x\geq m}\mu[x,+\infty)\int_m^x\frac{\omega'(u)^2}{h(u)}\,du<+\infty,
\end{equation}
where $m$ denotes the median of $\mu$. Moreover the optimal constant $C$ in (\ref{Poincaré}) denoted by $C_{\mathrm{opt}}$ verifies $$\max (D_\omega^-,D_\omega^+)\leq C_{\mathrm{opt}}\leq 4\max (D_\omega^-,D_\omega^+)$$
\end{prop}
This proposition follows at once from the celebrated Muckenhoupt criteria for the classical Poincar\'e inequality (see \cite{Muckenhoupt1972}). The following result completes the picture giving a large class of examples:
\begin{prop}\label{examplesdim1}
Let $\mu$ be an absolutely continuous probability measure on $\R$ with density $d\mu(x)=e^{-V(x)}\,dx$. Assume that the potential $V$ is of class $C^1$ and that $\omega$ verifies the following regularity condition:
\[\frac{\omega''(x)}{\omega'^2(x)}\xrightarrow[x\to +\infty]{} 0.\]
If $V$ is such that
\begin{equation}
\liminf_{x\to\pm\infty}\frac{\sgn(x)V'(x)}{\omega'(x)}>0,
\end{equation}
then the probability measure $\mu$ verifies the Poincar\'e inequality $\SG(\omega, C)$ for some $C>0$.
\end{prop}

In dimension $d$, one gets:
\begin{prop}\label{CSdimd}
Let $\mu$ be a probability measure on $\R^d$ absolutely continuous with respect to the Lebesgues measure, with $d\mu(x)=e^{-V(x)}\,dx$ with $V$ a function of class $C^2.$ Suppose that $\omega$ is of class $C^3$ on $\R$ and such that $\omega'(0)>0$ and
\[\forall x\in \R,\quad \left|\frac{\omega^{(3)}}{(\omega')^3}(x)\right|\leq M,\]
for some $M>0$. If there is some constant $u>0$ such that
\[\liminf_{|x|\to+\infty}\frac{1}{u^2}\sum_{i=1}^d\lc\frac{1}{4}\lp\frac{\partial V}{\partial x_i}\rp^2\lp\frac{x}{u}\rp-\frac{\partial^2V}{\partial x_i^2}\lp\frac{x}{u}\rp\rc\frac{1}{\omega'(x_i)^2}>dM,\]
then the probability measure $\mu$ satisfies $\SG(\tilde{\omega},C)$ for some $C$, where $\tilde{\omega}(x)=\omega(ux),$ for all $x\in \R$.
\end{prop}
This condition will be easily derived from the condition
$\liminf_{|x|\to+\infty}|\nabla V|(x)^2-\Delta V(x)>0,$
which is known to imply the classical Poincar\'e inequality.

\subsubsection{Links with Transportation-Cost inequalities} In Section \ref{Sec TCI}, we show the equivalence between the Poincar\'e inequalities for the metric $d_\omega$ and certain transportation-cost inequalities. Transportation-cost inequalities were first introduced by K. Marton and M. Talagrand in \cite{Mar86,Mar96} and \cite{Tal96a}. For recent advances in the understanding of these inequalities consult \cite{CG06}, \cite{GGM05}, \cite{Gozlan2007}, \cite{Wang04,Wang}. In these inequalities one tries to bound an optimal transportation cost in the sens of Kantorovich by the relative entropy functional. More precisely, if $c:\X\times\X\to\R^+$ is a measurable map on some metric space $\X$, the optimal transportation cost between $\nu$ and $\mu \in \PX$ (the set of probability measures on $\X$) is defined by
\[\mathcal{T}_c(\nu,\mu)=\inf_{\pi\in P(\nu,\mu)}\int c(x,y)\,d\pi,\]
where $P(\nu,\mu)$ is the set of probability measures $\pi$ on $\X\times\X$ such that $\pi(dx,\Y)=\nu(dx)$ and $\pi(\X,dy)=\mu(dy)$.
One says that $\mu$ satisfies the transportation cost inequality with the cost function $c(x,y)$ if
\begin{equation}\label{TCI}
\forall \nu\in \PX,\qquad \Tnm\leq \Hnm,
\end{equation}
where $\Hnm$ denotes the relative entropy of $\nu$ with respect to $\mu$ and is defined by $\Hnm=\int \log\lp\frac{d\nu}{d\mu}\rp\,d\nu$ if $\nu$ is absolutely continuous with respect to $\mu$ and $\Hnm=+\infty$ otherwise.

Transportation cost inequalities are known to have good tensorization properties and to yield concentration results independent of the dimension (all these facts are recalled in section \ref{Sec TCI}). For example, the celebrated $\T_2$ inequality which corresponds to cost functions of the form $(x,y)\mapsto a|x-y|_2^2$ gives gaussian concentration (see e.g \cite{Tal96a}). A celebrated result of Otto and Villani shows that the Lograithmic Sobolev inequality implies $\T_2$ (see \cite{OV00}).

Let us say that $\mu\in \mathcal{P}(\Rd)$ satisfies the inequality $\T(\omega,a)$ if it satisfies the transportation cost inequality (\ref{TCI}) with the cost function $(x,y)\mapsto\alpha\lp a\domxy\rp$

One proves the following
\begin{thm}\label{SG=TCI}
Let $\mu$ be a probability measure on $\Rd$ absolutely continuous with respect to Lebesgues measure with a positive density.
Then $\mu$ satisfies the Poincar\'e inequality $\SG(\omega,C)$ for some $C>0$ if and only if it satisfies the transportation-cost inequality $\T(\omega,a)$ for some $a>0$.\\
More precisely,
\begin{itemize}
\item[$\bullet$]if $\mu$ satisfies $\SG(\omega,C)$ then it satisfies  $\T(\omega,\frac{1}{\sqrt{C}\kappa}),$ with $\kappa=\sqrt{18e^{\sqrt{5}}}.$
\item[$\bullet$]if $\mu$ satisfies the inequality $\T(\omega,a)$, then $\mu$ satisfies the inequality $\SG(\omega,\frac{1}{2a^2})$.
\end{itemize}
\end{thm}
This theorem is an easy extension of a result by Bobkov, Gentil and Ledoux concerning the classical Poincar\'e inequality (see \cite[Corollary 5.1]{BGL01}). This extension is performed using a very simple contraction principle for transportation cost inequalities. The author previously used this technique in \cite{Gozlan2007} to characterize a large class of transportation cost inequalities on the real line.

\subsubsection{Comparison with Latala-Oleszkiewicz inequalities} In Section \ref{Sec Comp}, we compare the inequalities $\SG(\omega, C)$ to other functional inequalities including the ones introduced by R. Latala and K. Oleszkiewicz in \cite{Latala2000a}. Let $r\in[1,2]$, one says that a probability measure $\mu$ on $\Rd$ satisfies the inequality $\LO(r,C)$ if
\begin{equation}\label{LO}
\sup_{p\in (1,2)}\frac{\int f^2\,d\mu-\lp\int f^p\,d\mu\rp^{2/p}}{(2-p)^{2(1-1/r)}}\leq C \int |\nabla f|^2\,d\mu.
\end{equation}
It is well known that these inequalities interpolate between Poincar\'e and Log-Sobolev. For $r=1$, the inequality (\ref{LO}) is Poincar\'e inequality $\SG(C)$ and for $r=2$ it is equivalent to the Logarithmic-Sobolev inequality (see \cite[Corollary 1]{Latala2000a}). The $\LO(r,C)$ inequalities on $\R$ were completely characterized by Barthe and Roberto in \cite{Barthe2003}.

Recall that a probability measure $\mu$ on $\Rd$ verifies the Logarithmic-Sobolev inequality with constant $C$, if for all smooth $f$,
\begin{equation}\label{Log-Sob}
\ent_\mu(f^2)\leq C \int |\nabla f|^2\,d\mu
\end{equation}
where $\ent_\mu(f^2):=\int f^2\log f^2\,d\mu-\int f^2\,d\mu \log\lp\int f^2\,d\mu\rp$.

If $\mu$ verifies $\LO(r,C)$ then a concentration inequality of the same order as the one given in Proposition \ref{concentration omegap} holds (see \cite[Theorem 1]{Latala2000a}).
In fact, one has the following
\begin{thm}\label{LO>SG}
Let $r\in[1,2]$ ; if $\mu$ verifies the Latala-Oleszkiewicz inequality $\LO(r,C)$ for some $C>0$ then it satisfies the Poincar\'e inequality $\SG(\omega_r,\tilde{C})$ for some constant $\tilde{C}$, where $\omega_r(x)=x^r$ on $\R^+$.
\end{thm}
Moreover, a counter example of Cattiaux and Guillin shows that the Logarithmic-Sobolev inequality is strictly stronger than the inequality $\SG(\omega_2,C)$ (see Remark \ref{Cattiaux Guillin}).



\section{Weighted forms of the Poincar\'e inequality}
\label{Sec SG}
\subsection{Links with the classical Poincar\'e inequality}
\proof[Proof of Proposition \ref{equivalence}.]
Let us denote $|\nabla f|_\omega$ (resp. $|\nabla f|_2$) the length of the gradient computed with respect to the metric $\dom$ (see (\ref{length})).
If $f:\R^d\to\R$ is locally Lipschitz for the euclidean metric, then according to Rademacher theorem, one has
\[\limsup_{y\to x} \frac{|f(x)-f(y)|}{|x-y|_2}=\lc\sum_{i=1}^d \lp\frac{\partial f}{\partial x_i}\rp^2\rc^{1/2}=|\nabla f|_2(x),\]
for $\mu$ a.e. $x\in \R^d$, and so the length of the gradient equals the norm of the vector $\nabla f$ $\mu$ a.e.\\

Locally lipschitz function for $\dom$ and $|\,.\,|_2$ are related in the following way. A function $g:\R^d\to\R$ is locally Lipschitz for $d_\omega(\,.\,,\,.\,)$ if and only if $g\circ \omega^{-1}$ is locally Lipschitz for $|\,.\,|_2$.

[(i)$\Rightarrow$(ii)] Define $\tilde{\mu}=\omega^\sharp\mu$. Let $f:\R^d\to\R$ be locally Lipschitz for $|\,.\,|_2$, then $f\circ\omega$ is locally Lipschitz for $\dom$, and \[\Var_{\tilde{\mu}}(f)=\Var_\mu(f\circ\omega)\leq \int |\nabla (f\circ \omega)|^2_\omega\,d\mu\stackrel{(\ast)}{=}\int |\nabla f|^2_2\circ\omega=\int |\nabla f|_2^2\,d\mu,\]
where $(\ast)$ follows from the easy to check identity: $|\nabla (f\circ\omega)|_\omega=|\nabla f|_2\circ\omega$.

[(ii)$\Rightarrow$(i)] The proof is the same.

[(ii)$\Rightarrow$(iii)] Take $f:\R^d\to\R$ such that $f\circ\omega^{-1}$ is of class $C^1$.
Then
\[\Var_\mu(f)=\Var_{\tilde{\mu}}(f\circ\omega^{-1})\leq \int |\nabla (f\circ\omega^{-1})|_2^2\circ\omega\,d\mu=\int\sum_{i=1}^d \frac{1}{\omega'(x_i)^2}\lp\frac{\partial f}{\partial x_i}(x)\rp^2\,d\mu(x)\]

[(iii)$\Rightarrow$(ii)] Apply the weighted Poincar\'e inequality to the function $f\circ\omega$ with $f$ of class $C^1$.
\endproof
\subsection{Poincar\'e inequalities and concentration - the abstract case.}\label{sec concentration}
In order to recall how concentration estimates can be derived from the Poincar\'e inequality, let us briefly sketch the proof of Theorem \ref{BL1}.

\proof[Sketch of proof of Theorem \ref{BL1}]~

[First step] According to \cite[Theorem 3.1]{BL97} (which is the main result of \cite{BL97}), $\mu$ enjoys a modified Logarithmic-Sobolev inequality: for all $0<s<\frac{2}{\sqrt{C}}$ and for all locally Lipschitz $f:\X\to\R$ such that $|\nabla f|\leq s$ $\mu$ a.e. one has
\begin{equation}\label{modified logsob}
\ent_\mu(e^f)\leq L(s) \int |\nabla f|^2e^f\,d\mu,
\end{equation}
where $L(s)=\frac{C}{2}\lp\frac{2+\sqrt{C}s}{2-\sqrt{C}s}\rp^2e^{s\sqrt{5C}}.$

[Second step] \textit{Tensorization}. Thanks to the tensorization property of the entropy functional,
\[\ent_{\mu^n}(e^f)\leq \int \sum_{i=1}^n \ent_\mu(e^{f_i})\,d\mu^n,\]
for all $f:\X^n\to\R$.

Applying this inequality together with (\ref{modified logsob}) yields
\begin{equation}\label{modified logsob n}
\ent_{\mu^n}(e^f)\leq L(s) \int \sum_{i=1}^n |\nabla_i f|^2e^f\,d\mu,
\end{equation}
for all $0<s<\frac{2}{\sqrt{C}}$ and $f:\X^n\to\R$ such that $\max_{1\leq i\leq n}|\nabla_i f|\leq s$ $\mu^n$ a.e.

[Third step] \textit{Herbst argument.} Thanks to the homogeneity one can suppose that $f:\X^n\to\R$ is such that $\max_{1\leq i\leq n}|\nabla_i f|\leq 1$ ($b=1$) and $\sum_{i=1}^n |\nabla_if|^2\leq a^2$. Define $Z(\lambda)=\int e^{\lambda f}\,d\mu^n$. Then, applying (\ref{modified logsob n}) to $\lambda f$, one easily obtains the following differential inequality
\[\forall 0<\lambda\leq s<\frac{2}{\sqrt{C}},\qquad \frac{d}{d\lambda}\lp\frac{\log(Z(\lambda))}{\lambda}\rp \leq L(s)a^2,\]
and since $\frac{\log(Z(\lambda))}{\lambda}\xrightarrow[\lambda\to0]{} \int f\,d\mu^n$, one gets
\[\forall 0<\lambda\leq s<\frac{2}{\sqrt{C}},\qquad\int e^{\lambda f}\,d\mu^n\leq e^{\lambda^2 L(s)a^2+\lambda \int f\,d\mu^n}\]

[Fourth step] \textit{Chebischev argument.} This latter inequality on the Laplace transform yields via Chebischev argument:
\[\forall t\geq 0,\qquad \mu^n\lp f\geq\int f\,d\mu^n+t\rp\leq e^{-h_s(t)},\]
where \[h_s(t)=\sup_{\lambda\in [0,s]}\{\lambda t-L(s)a^2\lambda^2\}=\la\begin{array}{ll}\frac{t^2}{4L(s)a^2}&\text{if } 0\leq t\leq 2L(s)a^2s\\ st-L(s)a^2s^2 & \text{if } t\geq 2L(s)a^2s\end{array}\right.\]
Now it easy to see that,
$h_s(t)\geq \min\lp\frac{t^2}{4L(s)a^2},\frac{st}{2}\rp.$ For $s=1/\sqrt{C}$ one obtains after some computations,
\[h_s(t) \geq \min\lp \frac{t^2}{C\kappa^2a^2},\frac{t}{\sqrt{C}\kappa}\rp\qquad \text{with } \kappa=\sqrt{18e^{\sqrt{5}}}.\]
\endproof

\proof[Sketch of proof of Corollary \ref{Cor-concentration}.]
Take $A\subset \X^n$, such that $\mu^n(A)\geq 1/2$ and define $F(x)=\inf_{a\in A}\sum_{i=1}^n\alpha(d(x_i,a_i))$, where $\alpha(u)=\min(|u|,u^2)$.
Then for all $r\geq 0$, the function $f=\min(F,r)$ verifies (see the details in \cite{BL97}):
$\max_{1\leq i\leq n}|\nabla_i f|\leq 2$ and $\sum_{i=1}^n |\nabla_if|^2\leq 4r.$
Moreover since $\mu^n(A)\geq 1/2$, one has
$\int f\,d\mu^n = \int f \1_{A^c}\,d\mu^n\leq r(1-\mu^n(A))\leq r/2.$
Consequently, applying (\ref{deviation}) to $f$ yields:
\[\mu^n(F\geq r)=\mu^n(f\geq r)\leq \mu^n\lp f\geq \int f\,d\mu^n + r/2\rp\leq e^{-rK(C)},\]
with $K(C)=\frac{1}{16}\min\lp\frac{1}{C\kappa^2},\frac{1}{\sqrt{C}\kappa}\rp=\frac{1}{16}\alpha(\frac{1}{\sqrt{C}\kappa}).$
This achieves the proof of (\ref{concentration}).\endproof

\subsection{The $\SG(\omega,C)$ inequality and concentration.}
\begin{prop}\label{concentration omega}
Suppose that $\mu\in \mathcal{P}(\Rd)$ satisfies $\SG(\omega,C)$ for some $C>0$. Then for all $n\geq 1$ and all $A\subset \Rdn$, one has
\[\forall h\geq 0,\qquad \mu^n\lp A^h_\omega \rp\geq 1-e^{-K(C)h/d},\]
where $K(C)=\alpha\lp\frac{1}{\sqrt{C}\kappa}\rp/16$ and $A^h_\omega$ is defined by
\[A^h_\omega=\la (x_1,\ldots,x_n)\in \Rdn : \inf_{a\in A}\sum_{i=1}^n\sum_{j=1}^d \alpha\circ\omega\lp\frac{|x_{i,j}-a_{i,j}|}{2}\rp\leq h\ra.\]
(For all $1\leq i\leq n$, $x_{i,j}, 1\leq j\leq d$ are the coordinates of the vector $x_i\in \Rd.$)
\end{prop}
\begin{rem}
The fact that the dimension $d$ appears in the preceding result is not important. The important thing is that the constants do not depend on the dimension $n$.
\end{rem}

We need the following elementary lemmas:
\begin{lem}\label{super +}
If $f:\R^+\to\R$ is such that $x\mapsto f(x)/x$ is non decreasing then $f$ is super additive, that is to say: $f(x+y)\geq f(x)+f(y)$ for all $x,y\geq 0$.
\end{lem}
\proof
Let $0<x\leq y$ ; $f(x+y)=f(y(1+x/y))\geq (1+x/y)f(y)=f(y)+x f(y)/y\geq f(y)+xf(x)/x=f(y)+f(x).$
\endproof

\begin{lem}\label{elementary lemma}
For all $x,y\in \R$, $|\omega(x)-\omega(y)|\geq \omega\lp\frac{|x-y|}{2}\rp.$
\end{lem}
\proof
According to the Lemma \ref{super +}, the function $\omega$ is super additive on $\R^+$.
Let $x\geq y$.
If $x\geq y\geq 0$, then using the super additivity of $\omega$, one gets $\omega(x)=\omega((x-y)+y)\geq \omega(x-y)+\omega(y)$, so $\omega(x)-\omega(y)\geq \omega(x-y)\geq \omega((x-y)/2)$.
If $0\geq x\geq y$, then, according to the preceding case,  $\omega(x)-\omega(y)=\omega(-y)-\omega(-x)\geq \omega((-y+x)/2)=\omega((x-y)/2)$.
If $x\geq 0 \geq y$, then
$\omega(x)-\omega(y)=\omega(x)+\omega(-y)\geq \omega (\max(x,-y))\geq \omega ((x-y)/2)$.
\endproof

\begin{lem}\label{super *}
The function $\alpha(u)=\min(|u|,u^2)$ is such that $\alpha(au)\geq \alpha(a)\alpha(u)$, for all $a,u\geq 0$.
\end{lem}
\proof
If $0<a\leq 1$, then $\alpha(au)/a=u^2$ if $u\leq 1/a$ and $\alpha(au)/a=u/a$ if $u\geq 1/a$. If $u\leq 1$, one has $\alpha(au)/a=\alpha(u)$. If $u\in[1,1/a]$, then $u^2\geq u$ and so $\alpha(au)/a\geq \alpha(u)$. If $u\geq 1/a$, then $u/a\geq a$ and so $\alpha(au)/a\geq \alpha(u)$. The case $a\geq 1$ can be handled in a similar way.
\endproof
\proof[Proof of Proposition \ref{concentration omega}]
First, $\domxy\geq \frac{1}{\sqrt{d}}\sum_{i=1}^d |\omega(x_i)-\omega(y_i)|$, for all $x,y\in \R^d$.\\
Now,
\begin{align*}
\alpha(\domxy)&\geq \alpha\lp\sum_{i=1}^d \frac{1}{\sqrt{d}}|\omega(x_i)-\omega(y_i)|\rp \stackrel{(i)}{\geq} \sum_{i=1}^d \alpha\lp\frac{1}{\sqrt{d}}|\omega(x_i)-\omega(y_i)|\rp\\
&\stackrel{(ii)}{\geq} \sum_{i=1}^d \alpha\lp\frac{1}{\sqrt{d}}\omega\lp\frac{|x_i-y_i|}{2}\rp\rp
\stackrel{(iii)}{\geq} \frac{1}{d}\sum_{i=1}^d \alpha\circ\omega\lp\frac{|x_i-y_i|}{2}\rp
\end{align*}
where $(i)$ comes from the super additivity of the function $\alpha$, $(ii)$ from Lemma \ref{elementary lemma} and $(iii)$ from Lemma \ref{super *}.

Consequently,
\[\inf_{a\in A} \sum_{i=1}^n \alpha(d_\omega(x_i,a_i))\geq \frac{1}{d}\inf_{a\in A} \sum_{i=1}^n \alpha\circ\omega\lp\frac{|x_i-y_i|}{2}\rp.\]
Applying (\ref{concentration}) yields immediately the desired result.
\endproof

\proof[Proof of Proposition \ref{concentration omegap}.]
Suppose $p\in[1,2]$ ;
in view of Theorem \ref{concentration omega}, it is enough to prove that
\[\sum_{k=1}^{nd} \alpha\circ\omega_p(u_k)\leq h \Rightarrow u=(u_1,\ldots,u_{nd})\in \sqrt{h}B_2 + h^{1/p}B_p.\]
Let $v=(v_1,\ldots,v_{nd})$ and $w=(w_1,\ldots,w_{nd})$ be defined by $v_k=u_k$ if $u_k\in [-1,1]$ and $v_k=0$ if $|u_k|>1$ and $w=u-v$. Then,
\[\sum_{k=1}^{nd} \alpha\circ\omega_p(u_k)=|v|_2^2+|w|_p^p\leq h.\]
So, $|v|_2\leq \sqrt{h}$ and $|w|_p\leq h^{1/p}.$ Since $u=v+w$, one concludes that $u\in \sqrt{h}B_2 + h^{1/p}B_p.$

Now, if $p\geq 2$, then $\forall x\geq 0$, $\alpha\circ\omega(x)\geq x^2$ and $\forall x\geq 0$, $\alpha\circ\omega(x)\geq x^p$. This observation together with Theorem \ref{concentration omega} easily implies the result.
\endproof
\section{Workable sufficient conditions for $\SG(\omega,C)$.}\label{Sec CS}
\subsection{Dimension one.}
\proof[Proof of Proposition \ref{CSdim1}.]
According to Muckenhoupt criterion, a probability measure $d\nu=h\,dx$ having a positive continuous density with respect to Lebesgues measure, satisfies the classical Poincar\'e inequality if and only if
\[D^-=\sup_{x\leq m} \nu(-\infty,x]\int_x^m \frac{1}{h(u)}\,du<+\infty\quad\text{and}\quad D^+=\sup_{x\geq m} \nu[x,+\infty)\int_m^x \frac{1}{h(u)}\,du<+\infty,\]
and the optimal constant $C_{opt}$ verifies $\max(D^-,D^+)\leq C_{opt}\leq 4\max(D^-,D^+).$
Now, according to Proposition (\ref{equivalence}) $\mu$ satisfies $\SG(\omega,C)$ if and only if $\tilde{\mu}=\omega^\sharp\mu$ satisfies $\SG(C)$. The density of $\tilde{\mu}$ is $\tilde{h}=\frac{h\circ \omega^{-1}}{\omega'\circ \omega^{-1}}$. Plugging $\tilde{h}$ in Muckenhoupt conditions gives immediately the announced result.
\endproof

\proof[Proof of Proposition \ref{examplesdim1}.]
Let $\tilde{\mu}=\omega^\sharp\mu$ and let $\nu$ be the symmetric exponential probability measure on $\R$, that is the probability measure with density $d\nu(x)=\frac{1}{2}e^{-|x|}\,dx$. It is well known that it verifies the following Poincar\'e inequality:
\begin{equation}\label{Poincexp}
\Var_\nu(g)\leq 4 \int g'^2(x)\,d\nu(x),
\end{equation}
for all smooth $g$ (see for example \cite[Lemma 2.1]{BL97}).
Let $T:\R\to\R$ be the map defined by $T(x)=F_{\tilde{\mu}}^{-1}\circ F_\nu(x)$, with $F_\nu(x)=\nu(-\infty,x]$ and $F_{\tilde{\mu}}(x)=\tilde{\mu}(-\infty,x]$. It is well known that $T$ is increasing and transports $\nu$ on $\tilde{\mu}$ which means that $T^\sharp\nu=\tilde{\mu}$. Let us apply inequality (\ref{Poincexp}) to a function $g=f\circ T$. It yields immediately:
\[\Var_{\tilde{\mu}}(f)\leq 4\int f'^2\lp T'\circ T^{-1}\rp ^2\,d\tilde{\mu}
 \leq 4 \lp\sup_{x\in \R} T'(x)\rp^2\int f'^2\,d\tilde{\mu}\]
As a conclusion, if the map $T$ is $L$ Lipschtitz then $\tilde{\mu}$ verifies Poincar\'e inequality $\SG(4L^2)$. The probability $\tilde{\mu}$ has density $d\tilde{\mu}(x)=e^{-\tilde{V}(x)}\,dx$, with $\tilde{V}(x)=V(\omega^{-1}(x))+\log\omega'\circ\omega^{-1}(x).$ It is proved in \cite{Gozlan2007} (see Proposition 34) that a sufficient condition for $T$ to be Lipschitz is that $\liminf_{x\to\pm\infty} \sgn(x)\tilde{V}'(x)>0.$ But $\tilde{V}'(\omega(x))=\frac{V'(x)}{\omega'(x)}+ \frac{\omega''(x)}{\omega'^2(x)}$ and by assumption $\frac{\omega''(x)}{\omega'^2(x)}\to 0$ when $x$ goes to $\infty$. Thus $\liminf_{x\to\pm\infty} \sgn(x)\tilde{V}'(x)=\liminf_{x\to\pm\infty} \frac{\sgn(x)V'(x)}{\omega'(x)},$ which achieves the proof.
\endproof
\begin{rem}
The condition $\liminf_{x\to\pm\infty} \frac{\sgn(x)V'(x)}{\omega'(x)}>0$ can also be derived from Proposition \ref{CSdim1} using the same techniques as in e.g \cite[Theorem 6.4.3]{Log-Sob}. But this method has the disadvantage of introducing useless technical assumptions such as $\lim_{\pm\infty} V''/(V'^2)=0$.
\end{rem}
\begin{rem}\label{Cattiaux Guillin}
According to Theorem \ref{LO>SG}, the Logarithmic Sobolev inequality is stronger than the Poincar\'e inequality $\SG(\omega_2,C)$.
In \cite{CG06}, P. Cattiaux and A. Guillin were able to construct a potential $V$ on $\R$ satisfying $V(-x)=V(x)$ and $\liminf_{x\to +\infty}V'(x)/x>0$ but such that the probability measure $d\mu=e^{-V(x)}\,dx$ does not satisfy the Bobkov-Götze necessary and sufficient condition for the Logarithmic Sobolev inequality (see \cite{BG99}). According to Proposition \ref{examplesdim1}, this shows that the Logarithmic Sobolev inequality is strictly stronger than the inequality $\SG(\omega_2,C)$.
\end{rem}

\subsection{Dimension $d$.}
\proof[Proof of Proposition \ref{CSdimd}.]
It is well known that a probability $d\nu(x)=e^{-W(x)}\,dx$ on $\R^d$ satisfies the classical Poincar\'e inequality if $W$ verifies the following condition:
\begin{equation}\label{CSPoinc}
\liminf_{|x|\to +\infty}\ |\nabla W|^2(x)-\Delta W(x)>0
\end{equation}
Suppose that $\mu$ is an absolutely continuous probability measure on $\R^d$ with density $d\mu(x)=e^{-V(x)}\,dx$ with $V$ of class $C^2$. Then $\tilde{\mu}=\omega^\sharp\mu$ has density $d\tilde{\mu}(x)=e^{-\tilde{V}(x)}\,dx,$ with \[\forall x\in \R^d,\qquad\tilde{V}(x)=V(\omega^{-1}(x))+\sum_{i=1}^d \log \omega'\circ\omega^{-1}(x_i).\]
According to Proposition \ref{equivalence}, to show that $\mu$ satisfies the inequality $\SG(\omega,C)$ for some $C>0$ it is enough to show that $\tilde{\mu}$ satisfies the inequality $\SG(C)$ and a sufficient condition for this is that $\tilde{V}$ fulfills condition (\ref{CSPoinc}).

Elementary computations yield
\begin{align*}
\frac{\partial\tilde{V}}{\partial x_i}(\omega(x))&=\frac{1}{\omega'(x_i)}\frac{\partial V}{\partial x_i}(x)+\frac{\omega''(x_i)}{\omega'^2(x_i)}\\
\frac{\partial^2\tilde{V}}{\partial x_i^2}(\omega(x))&=-\frac{\omega''(x_i)}{\omega'^3(x_i)}\frac{\partial V}{\partial x_i}(x)+\frac{1}{\omega'^2(x_i)}\frac{\partial^2V}{\partial x_i^2}(x)+\frac{\omega^{(3)}(x_i)}{\omega'^3(x_i)}-2\frac{\omega''^2(x_i)}{\omega'^4(x_i)}
\end{align*}

Let $I(x)=|\nabla \tilde{V}|^2(\omega(x))-\Delta \tilde{V} (\omega(x))$ ; one has:
\[I(x)=\sum_{i=1}^d \frac{1}{\omega'^2(x_i)}\lc\lp\frac{\partial V}{\partial x_i}\rp^2(x)-\frac{\partial^2V}{\partial x_i^2}(x)\rc +3 \sum_{i=1}^d \frac{\omega''(x_i)}{\omega'^3(x_i)}\frac{\partial V}{\partial x_i}(x)+3\sum_{i=1}^d\frac{\omega''^2(x_i)}{\omega'^4(x_i)}-\sum_{i=1}^d\frac{\omega^{(3)}(x_i)}{\omega'^3(x_i)}.\]
Using the inequality $uv\geq -u^2-v^2/4$, one has
\begin{align*}
3 \sum_{i=1}^d \frac{\omega''(x_i)}{\omega'^3(x_i)}\frac{\partial V}{\partial x_i}(x)&=3\sum_{i=1}^d \lp \frac{\omega''(x_i)}{\omega'^2(x_i)}\rp\cdot\lp\frac{1}{\omega'(x_i)}\frac{\partial V}{\partial x_i}(x)\rp\\
&\geq -3\sum_{i=1}^d\frac{\omega''^2(x_i)}{\omega'^4(x_i)}-\frac{3}{4}\sum_{i=1}^d \frac{1}{\omega'^2(x_i)}\lp\frac{\partial V}{\partial x_i}\rp^2(x),
\end{align*}
and so
\[I(x)\geq \sum_{i=1}^d \frac{1}{\omega'^2(x_i)}\lc\frac{1}{4}\lp\frac{\partial V}{\partial x_i}\rp^2(x)-\frac{\partial^2V}{\partial x_i^2}(x)\rc-\sum_{i=1}^d \frac{\omega^{(3)}(x_i)}{\omega'^3(x_i)}.\]
Since, $\liminf_{|x|\to+\infty} I(x)=\liminf_{y\to+\infty}|\nabla \tilde{V}|^2(y)-\Delta \tilde{V} (y)$ and $\sum_{i=1}^d \frac{\omega^{(3)}(x_i)}{\omega'^3(x_i)}\leq dM$, one concludes that $\tilde{V}$ satisfies (\ref{CSPoinc}) as soon as
\[\liminf_{|x|\to +\infty}\sum_{i=1}^d \frac{1}{\omega'^2(x_i)}\lc\frac{1}{4}\lp\frac{\partial V}{\partial x_i}\rp^2(x)-\frac{\partial^2V}{\partial x_i^2}(x)\rc>dM.\]
Applying this latter condition to the probability measure $\mu_u=(u\Id)^\sharp\mu$, (where $\Id$ is the identity function) which has density $d\mu_u(x)=\frac{1}{u^d}e^{-V(x/u)}\,dx$ gives the condition of Proposition \ref{CSdimd}.
\endproof
\section{Transportation-cost inequalities}\label{Sec TCI}
\subsection{Basic properties}
\begin{prop}[Tensorization]\label{tensorization TCI}
Suppose that $\mu\in\mathcal{P}(\X)$ satisfies the transportation cost inequality (\ref{TCI}) with the cost function $c(x,y)$, then $\mu^n$ satisfies the transportation cost inequality on $\X^n$ with the cost function $c^{\oplus n}(x,y)=\sum_{i=1}^n c(x_i,y_i)$.
In other words,
\[\forall \nu\in\mathcal{P}(\X^n),\quad \inf_{\pi \in P(\nu,\mu^n)}\int \sum_{i=1}^n c(x_i,y_i)\,d\pi \leq \operatorname{H}(\nu \mid\mu^n),\]
where $P(\nu,\mu^n)$ is the set of probability measures on $\X^n\times\X^n$ such that $\pi(dx,\X^n)=\nu(dx)$ and $\pi(\X^n,dy)=\mu^n(dy)$.
\end{prop}
This result goes back to the first works of K. Marton on the subject (see \cite{Mar86,Mar96}). A proof can be found in \cite{Gozlan2007a}.

Let us explain how to derive concentration inequalities from the inequality $\T(\omega,a)$.
\begin{prop}\label{concentration TCI}
If $\mu$ satisfies the transportation cost inequality $\T(\omega,a)$, then for all $n\geq 1$ and all $A\subset \R^{nd}$,
\[\forall h\geq 0,\qquad \mu^{n}\lp A^h_\omega\rp\geq 1-\frac{1}{\mu^n(A)}e^{-h\alpha(a/\sqrt{d})/2},\]
where the enlargement is defined by
\[A^h_\omega=\la y=(y_1,\ldots,y_n)\in \Rdn  : \inf_{x\in A}\sum_{i=1}^n \sum_{j=1}^n \alpha\circ\omega \lp\frac{x_{i,j} - y_{i,j}}{2}\rp\leq h \ra.\]
\end{prop}
\begin{rem}
According to Theorem \ref{SG=TCI}, if $\mu$ satisfies the inequality $\SG(\omega,C)$ then it satisfies $\T(\omega, a)$ with $a=\frac{1}{\sqrt{C}\kappa}$. With this value of $a$ the concentration inequality given by Proposition \ref{concentration TCI} is almost the same as the one derived in Proposition \ref{concentration omega}.
\end{rem}
We will need the following lemma:
\begin{lem}\label{elementary lemma2}
The function $\alpha(u)=\min(|u|,u^2)$ is such that
 $\alpha(x+y)\leq 2(\alpha(x)+\alpha(y))$, for all $x,y\geq 0.$
\end{lem}
\proof
If $x+y\leq 1$, then $\alpha(x+y)=(x+y)^2\leq 2(x^2+y^2)=2(\alpha(x)+\alpha(y))$.\\
Now, suppose that $x+y\geq 1$. \\
If $x\leq 1$ and $y\leq 1$, then $\alpha(x+y)=x+y\leq (x+y)^2\leq 2(x^2+y^2)=2(\alpha(x)+\alpha(y)).$\\
If $x\leq 1$ and $y\geq 1$, then
$x\leq y \Rightarrow x-2x^2\leq y \Rightarrow x+y\leq 2(x^2+y)\Rightarrow \alpha(x+y)\leq 2(\alpha(x)+\alpha(y)).$
If $x\geq 1$ and $y\geq 1$, then $\alpha(x+y)=x+y=\alpha(x)+\alpha(y)\leq 2(\alpha(x)+\alpha(y)).$
\endproof
\proof[Proof of Proposition \ref{concentration TCI}.]
If $\mu$ satisfies $\T(\omega,a)$ on $\R^d$ then according to Theorem  (\ref{tensorization TCI}), $\mu^n$ satisfies the transportation cost inequality on $\Rdn$ with the cost function $c$ defined by
\[c :\lp(x_1,\ldots,x_n),(y_1,\ldots,y_n)\rp\in\Rdn \times \Rdn  \mapsto \sum_{i=1}^n \alpha(ad_\omega(x_i,y_i)).\]
Using the triangle inequality for the metric $\dom$ and Lemma \ref{elementary lemma2}, one has
\[\forall x,y,z \in \Rdn ,\qquad c(x,z)\leq  2c(x,y) + 2c(y,z).\]

Now, let $\nu_1$ and $\nu_2$ be two probability measures on $\R^{nd}.$ Take $\pi_1\in P(\nu_1,\mu^n)$ and $\pi_2\in P(\mu^n,\nu_2)$, then one can construct three random variables $X,Y,Z$ such that
$\mathcal{L}(X,Y)=\pi_1$ and $\mathcal{L}(Y,Z)=\pi_2$ (see for instance the Gluing Lemma of \cite{Vill} p. 208).
Then, one has
\begin{align*}
\mathcal{T}_{ c }(\nu_1,\nu_2)&\leq \E\lc c(X,Z)\rc\leq 2\E\lc c(X,Y)\rc + 2\E\lc c(Y,Z)\rc\\
&=2\int c(x,y)\,d\pi_1(x,y)+2\int c(y,z)\,d\pi_2(y,z).
\end{align*}
Optimizing on $\pi_1$ and $\pi_2$ gives
\[\mathcal{T}_{c}(\nu_1,\nu_2)\leq 2\mathcal{T}_{c}(\nu_1,\mu^n)+2\mathcal{T}_{c}(\nu_2,\mu^n)\]
Consequently, $\mu^n$ satisfies the following symmetrized transportation cost inequality:
\[\forall \nu_1,\nu_2 \in \mathcal{P}(\R^{nd}),\qquad \mathcal{T}_{c}(\nu_1,\nu_2)\leq 2\operatorname{H}(\nu_1\mid \mu^n)+2\operatorname{H}(\nu_2\mid \mu^n).\]
Take $d\nu_1=\1_A d\mu^n$ and $d\nu_2=\1_B\,d\mu^n$, then
\begin{align*}
\inf_{x\in A, y\in B} c(x,y)&\leq \mathcal{T}_{c}(\nu_1,\nu_2)\leq 2\operatorname{H}(\nu_1\mid \mu^n)+2\operatorname{H}(\nu_2\mid \mu^n)\\&=2\log(1/\mu^n(A))+2\log(1/\mu^n(B))
\end{align*}
Letting $c(A,B)=\inf_{x\in A, y\in B} c(x,y)$, one gets
\[\mu^{(n)}(A)\mu^{(n)}(B)\leq e^{-c(A,B)/2}.\]
Defining
$B=\{y : \inf_{x\in A} c(x,y)> h \}$
one gets $\mu^n(B)\leq \frac{1}{\mu^n(A)}e^{-h/2}$.
To obtain the announced inequality it is thus enough to compare  $A_\omega^h$ and $B$.
Take $x=(x_1,\ldots,x_n)\in \Rdn $ and $y=(y_1,\ldots,y_n)\in \Rdn $ ; then for all $i\in{1,\ldots,n}$, one has
\begin{align*}
\alpha\lp ad_\omega(x_i,y_i)\rp&\stackrel{(a)}{\geq}\alpha\lp\frac{a}{\sqrt{d}}\sum_{j=1}^d|\omega(x_{i,j})-\omega(y_{i,j})|\rp\stackrel{(b)}{\geq} \sum_{j=1}^d\alpha\lp\frac{a}{\sqrt{d}}|\omega(x_{i,j})-\omega(y_{i,j})|\rp\\
&\stackrel{(c)}{\geq} \sum_{j=1}^d\alpha\lp\frac{a}{\sqrt{d}}\omega\lp \frac{x_{i,j}-y_{i,j}}{2}\rp\rp \stackrel{(d)}{\geq} \alpha\lp a/\sqrt{d}\rp\sum_{j=1}^d\alpha\circ\omega\lp \frac{x_{i,j}-y_{i,j}}{2}\rp
\end{align*}
where
(a) follows from the comparison between the norms $|\,.\,|_2$ and $|\,.\,|_1$ in $\R^d$, (b) from Lemma \ref{super +}, (c) from Lemma \ref{elementary lemma} and (d) from Lemma \ref{super *}.

Consequently, if $\inf_{x\in A} \sum_{i=1}^n\sum_{j=1}^d\alpha\circ\omega\lp \frac{x_{i,j}-y_{i,j}}{2\sqrt{d}}\rp\geq h/\alpha(a/\sqrt{d})$, then $y$ belongs to $B$. From this follows that $\mu^n(A_\omega^h)\geq 1-\frac{1}{\mu^n(A)}e^{-\alpha(a/\sqrt{d})h/2},$ which achieves the proof.
\endproof
\begin{rem}
The idea of deriving concentration estimates from transportation cost inequalities goes back to Marton seminal work \cite{Mar86}. The above proof is essentially due to Talagrand (see the proof of \cite[Corollary 1.3]{Tal96a}).
\end{rem}
\subsection{Links with Poincar\'e inequality}
The proof of Theorem \ref{SG=TCI} relies on two ingredients. The first one is the following result by Bobkov, Gentil and Ledoux:
\begin{thm}[Bobkov,Gentil, Ledoux]\label{BGL} If a probability measure $\mu$ on $\R^d$ satisfies $\SG(C)$ then it satisfies the transportation cost inequality for the cost function $(x,y)\mapsto\alpha_s(|x-y|_2)$ for all $s<\frac{2}{\sqrt{C}}$, where
\[\alpha_s(t)=\la\begin{array}{ll}\frac{t^2}{4L(s)}&\text{if } |t|\leq 2L(s)s\\
s|t|-L(s)s^2 & \text{otherwise} \end{array}\right.\qquad{with}\quad L(s)=\frac{C}{2}\lp\frac{2+\sqrt{C}s}{2-\sqrt{C}s}\rp^2e^{s\sqrt{5C}}.\]
\end{thm}
In particular, if one takes $s=\frac{1}{\sqrt{C}}$,
then it is easy to check that $\alpha_s(t)\geq \alpha\lp \frac{t}{\sqrt{C}\kappa}\rp$, where $\alpha(u)=\min(|u|,u^2)$ and $\kappa=\sqrt{18e^{\sqrt{5}}}$.
Thus if $\mu$ satisfies $\SG(C)$ it satisfies the transportation cost inequality with the cost function $(x,y)\mapsto\alpha\lp\frac{|x-y|_2}{\sqrt{C}\kappa}\rp.$
In other words, with the definition of the transportation cost inequality $\T\lp\omega,a\rp$, the preceding result can be restated as follows
\begin{cor}
If $\mu$ is a probability measure on $\R^d$ satisfying the classical Poincar\'e inequality $\SG(C)$ for some $C>0$, then it satisfies the transportation-cost inequality $\T\lp \Id,\frac{1}{\sqrt{C}\kappa}\rp.$ (where $\Id:\R\to\R : x\mapsto x$ is the identity function.)
\end{cor}
The converse is also true:
\begin{prop}\label{TCI->Poinc}
If $\mu$ satisfies $\T\lp \Id,a\rp,$ for some $a>0$, then $\mu$ satisfies the inequality $\SG(\frac{1}{2a^2})$.
\end{prop}
The proof of Proposition \ref{TCI->Poinc} is classical and can be found in various places (see e.g the proofs of \cite[Corollary 5.1]{BGL01} or \cite[Corollary 3]{Mau91}).

The second argument is a very simple contraction principle:
\begin{prop}
Let $\mu$ be a probability measure on a metric space $\X$ ; if $\mu$ satisfies the transportation cost inequality with the cost function $c:\X\times\X \to\R^+$, and if $T:\X\to\Y$ is a measurable bijection then, $T^\sharp\mu$ satisfies the transportation cost inequality with the cost function $(x,y)\mapsto c(T^{-1}(x),T^{-1}(y))$.
\end{prop}
This contraction principle goes back to Maurey's work on infimum convolution inequalities (see \cite{Mau91}). A proof can also be found in \cite{Gozlan2007}, where this simple property was intensively used to derive necessary and sufficient conditions for transportation cost inequalities on the real line.

Now let us apply the contraction principle together with Theorem
\ref{BGL} to prove that Poincar\'e inequalities $\SG(\omega,C)$ and transportation-cost inequalities $\T(\omega,a)$ are qualitatively equivalent.
\proof[Proof of Theorem \ref{SG=TCI}.]
If $\mu$ satisfies $\SG(\omega,C)$, then according to Proposition \ref{equivalence}, $\omega^\sharp\mu$ satisfies the classical Poincar\'e inequality $\SG(C)$, and according to Theorem \ref{BGL}, this implies that $\omega^\sharp\mu$ satisfies $\T(\Id,a)$, with $a=\frac{1}{\sqrt{C}\kappa}$. According to the contraction principle, $\mu$ (which is the image of $\omega^\sharp\mu$ under the map $\omega^{-1}$) satisfies the transportation cost inequality with the cost function $(x,y)\mapsto \alpha\lp a|\omega(x)-\omega(y)|_2\rp=\alpha\lp a\domxy\rp$ by definition of the metric $\dom$ (see (\ref{metric})).

Now suppose that $\mu$ satisfies $\T(\omega,a)$ for some $a>0.$ According to the contraction principle, $\omega^\sharp\mu$ satisfies $\T(\Id,a)$, and according to Proposition \ref{TCI->Poinc}, this implies that $\omega^\sharp\mu$ satisfies $\SG(\frac{1}{2a^2})$. Using Proposition \ref{equivalence}, one concludes that $\mu$ satisfies $\SG(\omega,\frac{1}{2a^2}).$ This achieves the proof.
\endproof

\begin{rem}
If $\mu$ satisfies the inequality $\T(\omega,a)$, it is easy to show that it verifies the transportation cost inequality with the cost function \[\Rd\times\Rd\to\R^+:(x,y)\mapsto \alpha(a/\sqrt{d})\sum_{i=1}^d \alpha\circ\omega\lp\frac{x_i-y_i}{2}\rp.\]
In particular, the inequality $\SG(\omega_2,C)$ implies Talagrand's $\T_2$ inequality, that is to say the transportation cost inequality with a cost function of the form $(x,y)\mapsto a|x-y|_2^2$ for some $a>0$. We do not know if the converse is true.
\end{rem}

\section{Comparison with other functional inequalities}
\label{Sec Comp}
In this section we will perform a comparison between the inequalities $\SG(\omega,C)$ and generalized Beckner-Latala-Oleszkiewicz inequalities introduced in \cite{Wang2005} and \cite{Barthe2006}.
\begin{defi}
Let $T:[0,1]\to\R^+$ be a non decreasing function and $\mu$ be a probability measure on $\R^d$. One says that $\mu$ satisfies the generalized Beckner-Latala-Oleszkiewicz inequality with the function $T$ and the constant $C>0$, if for all smooth $f$, one has
\begin{equation}\label{BLO}
\sup_{p\in (1,2)}\frac{\int f^2\,d\mu-\lp\int |f|^{p}\,d\mu\rp^{2/p}}{T(2-p)}\leq C \int |\nabla f|^2\,d\mu.
\end{equation}
If $\mu$ verifies (\ref{BLO}) one will say for short that $\mu$ satisfies the inequality $\BLO(T,C)$.
\end{defi}
The $\LO(r,C)$ inequality corresponds to the function $T(u)=u^{2(1-1/r)}.$

Dimension free concentration results can be deduced from the inequality $\BLO(T,C)$. The following result follows easily from Proposition 29 and Corollary 30 of \cite{Barthe2006}.
\begin{thm}\label{concentration BLO}
Let $T:[0,1]\to \R^+$ be a non decreasing function. Define $T(x)=T(1)$ for all $x\geq 1$ and let $\omega_T:\R\to \R$ be such that $\omega_T(-x)=-\omega(x)$ for all $x\in \R$ and
\begin{equation}\label{omega_T}
\forall t\geq 0,\quad \omega_T^{-1}(t)=\int_0^t \sqrt{T}(1/u)\,du.
\end{equation}
If $\mu$ satisfies the inequality $\BLO(T,C)$, then for all $n\geq 1$ and for all 1-Lipschitz function $f$:
\[\forall t\in \R^+\setminus \lc\sqrt{T(1)},2\sqrt{T(1)}\rc,\quad \mu^n\lp f\geq \int f\,d\mu+r\rp\leq e^{-\alpha\circ\omega_T(t/(3\sqrt{C}))}.\]
\end{thm}

We are going to prove the following result:
\begin{thm}\label{BLO>SG}
Let $T:[0,1]\to\R^+$ be a non-decreasing function such that $x\mapsto T(x)/x$ is non-increasing. If the measure $\mu$ verifies the inequality $\BLO(T,C)$ for some constant $C$ then it satisfies the inequality $\SG(\omega_T,\tilde{C})$.
\end{thm}

Let us admit Theorem \ref{BLO>SG} and let us prove Theorem \ref{LO>SG}.
\proof[Proof of Theorem \ref{LO>SG}.]
As noticed above, the inequality $\LO(r,C)$ is the same as $\BLO(T,C)$ with $T(u)=u^{2(1-1/r)}$.
According to Theorem \ref{BLO>SG}, $\mu$ verifies the inequality $\SG(\omega_T,C)$ for some $C$, where $\omega_T$ is given by (\ref{omega_T}).
A simple computation gives
$\omega_T(t)=t$ if $t\in[0,1]$ and $\omega_T(t)=t^r/r+1-1/r$, if $t\geq 1$. Thus, $\omega'_T(t)=\max(1,t^{r-1})$. On the other hand, $\omega_r'(t)=1$, if $t\in [0,1]$ and $\omega_r'(t)=rt^{r-1}$. Thus, $\frac{1}{r}\omega_r'(t)\leq \omega'_T(t)\leq \omega'_r(t)$, for all $t\geq 0$. Using (\ref{weightedPoincaré}), one concludes that
$\mu$ verifies $\SG(\omega_T,C)$ for some $C$ if and only if $\mu$ verifies $\SG(\omega_r,\tilde{C})$ for some $\tilde{C}$. This achieves the proof.
\endproof
The proof of this theorem relies on the capacity-measure formulation of the generalized Beckner-Latala-Oleszkiewicz inequalities due to Barthe, Cattiaux and Roberto \cite{Barthe2006}.

Let us recall the definition of a capacity-measure inequality (a good reference for this type of inequalities is the book of Maz'ja \cite{Maz'ja1985}).
\begin{defi}
Let $\mu$ be a probability measure on $\R^d$. Let $A\subset\Omega$ be Borel sets. One defines
\[\Capa_\mu(A,\Omega)=\inf \la\int |\nabla f|^2\,d\mu ; \1_A \leq f\leq \1_\Omega\ra.\]
The capacity of a set $A$ with $\mu(A)\leq 1/2$ is defined by
\begin{align*}
\Capa_\mu(A)&=\inf \la \Capa_\mu(A,\Omega) : A\subset \Omega \text{ and } \mu(\Omega)\leq 1/2\ra\\
&=\inf \la \int |\nabla f|^2\,d\mu ; f:\R^d\to[0,1],\ f_{|A}=1\text{ and } \mu(f=0)\geq 1/2\ra
\end{align*}
One says that $\mu$ satisfies a capacity-measure inequality if there is a function $\Theta:[0,1]\to \R^+$ and a constant $D>0$ such that for all $A$ with $\mu(A)\leq 1/2$,
\[\Theta(\mu(A))\leq D\Capa_\mu(A).\]
\end{defi}
The following theorem due to Barthe, Cattiaux and Roberto (see Theorem 18 and Lemma 19 of \cite{Barthe2006}) gives a capacity-measure transcription of the inequality $\BLO(T,C)$.
\begin{thm}\label{BLO<>Cap}
Let $T:[0,1]\to\R^+$ be a non-decreasing function such that $x\mapsto T(x)/x$ is non-increasing. Let $C>0$ be the optimal constant such that $\mu$ verifies the inequality $\BLO(T,C)$. Then $1/6 D\leq C\leq 20 D$, where $D$ is the optimal constant such that for all $A\subset \R^d$ with $\mu(A)\leq 1/2$, one has
\[\Theta(\mu(A))\leq D\Capa_\mu(A),\]
where $\Theta:\R^+\to \R^+$ is defined by:
\begin{equation}\label{Theta}
\forall x\in \R^+,\quad \Theta(x)=x\frac{1}{T\lp\frac{1}{\log\lp1+\frac{1}{x}\rp}\rp},
\end{equation}
with the convention that $T(x)=T(1)$ for $x\geq 1$.
\end{thm}
\begin{rem}
In fact we will only use the fact that the inequality $\BLO(T,C)$ implies the measure-capacity inequality $\Theta(\mu(A))\leq 6C \Capa_\mu(A)$, for all $A$ such that $\mu(A)\leq 1/2$. This is the easiest part of Theorem \ref{BLO<>Cap}.
\end{rem}
To prove Theorem \ref{BLO>SG}, one needs the following basic properties:
\begin{lem}If $T:[0,1]\to\R^+$ is a non-decreasing function such that $x\mapsto T(x)/x$ is non-increasing then the function $\Theta$ defined by (\ref{Theta}) is non-decreasing and verifies $\Theta(x+y)\leq \Theta(x)+\Theta(y)$ for all $x,y \in \R^+$.
\end{lem}
\proof
Let us write:
\[\Theta(x)=\frac{x}{h(x)}\cdot \frac{h(x)}{T(h(x))}\quad \text{with } h(x)=\frac{1}{\log(1+1/x)}.\]
The function $h$ is non-decreasing, and since $u\mapsto \frac{u}{T(u)}$ is non-decreasing, one concludes that $x\mapsto \frac{h(x)}{T(h(x))}$ is non decreasing. On the other hand, it is easy to see that the function $x\mapsto \frac{x}{h(x)}=x\log(1+1/x)$ is non-decreasing. As a product of non-decreasing and non-negative functions, the function $\Theta$ is itself non-decreasing.

Take $x\geq y>0$ ; using the fact that the function $x\mapsto \Theta(x)/x$ is non-increasing, one gets
\[\Theta(x+y)=\Theta(x(1+y/x))\leq (1+y/x)\Theta(x)=\Theta(x)+y \Theta(x)/x\leq \Theta(x)+\Theta(y).\]
This achieves the proof.
\endproof
Another ingredient of the proof is the following lemma which explains how behave capacity-measure inequalities under push-forward:
\begin{lem}\label{Cap-tilde}
Suppose that $\mu$ satisfies the capacity-measure inequality
\[\forall A\text{ with } \mu(A)\leq 1/2,\quad \Psi(\mu(A))\leq D\Capa_\mu(A).\]
Then $\tilde{\mu}=\omega^\sharp \mu$ verifies the inequality
\[\forall A\text{ with } \tilde{\mu}(A)\leq 1/2,\quad \Psi(\tilde{\mu}(A))\leq D\overline{\Capa}_{\tilde{\mu}}(A),\]
where \[\overline{\Capa}_{\tilde{\mu}}=\inf \la\int \sum_{i=1}^d \lp\omega'\circ\omega^{-1}(x_i)\rp^2\lp\frac{\partial f}{\partial x_i}\rp^2(x)\,d\tilde{\mu} ; f:\R^d \to [0,1],\ f_{|A}=1 \text{ and } \tilde{\mu}(f=0)\geq 1/2\ra.\]
\end{lem}
\proof
Let $A$ be such that $\tilde{\mu}(A)\leq 1/2$, and $f$ be such that $f=1$ on $A$ and $\tilde{\mu}(f=0)\geq 1/2$. Define $B=\omega^{-1}(A)$ and $g=f\circ\omega$. Then $\mu(B)=\tilde{\mu}(A)\leq 1/2$, $g\geq 1$ on $B$ and $\la g=0\ra =\omega^{-1}\la f=0\ra$ and so $\mu(g=0)=\tilde{\mu}(g=0)\geq 1/2$.
Applying the capacity-measure inequality verified by $\mu$ to $B$ and $g$ yields
\[\tilde{\mu}(A)=\mu(B)\leq D \int |\nabla g|^2\,d\mu=D\int \sum_{i=1}^d \lp\omega'\circ\omega^{-1}(x_i)\rp^2\lp\frac{\partial f}{\partial x_i}\rp^2(x)\,d\tilde{\mu}.\]
Optimizing over such functions $f$ gives the announced inequality for $\tilde{\mu}.$
\endproof
The next lemma explains how to compare the capacity $\overline{\Capa}_{\tilde{\mu}}$ to the usual capacity $\Capa_\mu$:
\begin{lem}\label{comp-Cap}
Let $B_\infty(r)=\la x\in \R^d : \max_{1\leq i\leq d} |x_i|\leq r\ra$, for all $r\geq 0$. If $A\subset B_\infty(r)$ and $\mu(A)\leq 1/2$, then
\[\overline{\Capa}_{\tilde{\mu}}(A)\leq 2\lp\omega'\circ \omega^{-1}(r+1)\rp^2\Capa_{\tilde{\mu}}(A)+8d\tilde{\mu}(B_\infty(r)^c).\]
\end{lem}
\proof
Let \[\Capa_{\tilde{\mu}}^r(A)=\inf \la \int |\nabla f|^2\,d\tilde{\mu} ; \1_A \leq f \leq \1_{B_\infty(r+1)} \text{ and } \tilde{\mu} (f=0)\geq 1/2 \ra.\]
Using the fact that the function $\omega'\circ\omega^{-1}$ is increasing on $\R^+$, one clearly has:
\[\overline{\Capa}_{\tilde{\mu}}(A)\leq \lp\omega'\circ \omega^{-1}(r+1)\rp^2\Capa_{\tilde{\mu}}^r(A).\]
Now let $f:\R^d\to [0,1]$ be such that $f_{|A}=1$ and $\tilde{\mu}(f=0)\geq 1/2.$ One can easily construct a cut-off function $\varphi$ such that: $\1_{B_\infty(r)}\leq \varphi\leq \1_{B_\infty(r+1)}$, and such that $\left|\frac{\partial \varphi}{\partial x_i}(x)\right|\leq 2$ for all $x\in \R^d$ and all $i$.
Let $g=f\varphi$; one has $\1_A\leq g\leq \1_{B_\infty(r+1)}$, $\tilde{\mu}(g=0)\geq \tilde{\mu}(f=0)\geq 1/2$ and
\begin{align*}
\Capa_{\tilde{\mu}}^r(A)&\leq \int |\nabla g|^2\,d\tilde{\mu} =\int |\nabla f\varphi + f\nabla \varphi|^2\,d\tilde{\mu}\\
&\leq  2\int |\nabla f|^2\varphi^2\,d\tilde{\mu} + 2\int f^2|\nabla \varphi|^2\,d\tilde{\mu}\\
&\leq 2\int |\nabla f|^2\,d\tilde{\mu} + 8d\tilde{\mu}(B_\infty(r)^c).\\
\end{align*}
Optimizing over $f$ yields:
\[\Capa_{\tilde{\mu}}^r(A)\leq 2 \Capa_{\tilde{\mu}}+8d\tilde{\mu}(B_\infty(r)^c).\]
\endproof

\proof[Proof of Theorem \ref{BLO>SG}.]
Define $\tilde{\mu}=\omega_T^\sharp \mu$. One wants to prove that $\tilde{\mu}$ verifies the classical Poincar\'e inequality.
According to Theorem \ref{BLO<>Cap}, the probability measure $\mu$ satisfies the capacity-measure inequality
\begin{equation}\label{eq-Cap-meas}
\forall A \text{ with } \mu(A)\leq 1/2,\quad \Theta(\mu(A))\leq 6C\Capa_\mu(A).
\end{equation}
According to Lemma \ref{Cap-tilde}, $\tilde{\mu}$ satisfies the capacity-measure type inequality:
\[\forall A \text{ with } \tilde{\mu}(A)\leq 1/2,\quad \Theta(\tilde{\mu}(A))\leq 6C\overline{\Capa}_{\tilde{\mu}}(A),\]
where $\overline{\Capa}_{\tilde{\mu}}$ is defined in the lemma.

Let $B_\infty(r)=\la x\in \R^d : \max_{1\leq i\leq d}(|x_i|)\leq r \ra$, for all $r\geq 0$.
Let $A\subset \R^d$ with $\tilde{\mu}(A)\leq 1/2$; one has
\begin{align*}
\Theta(\tilde{\mu}(A))&\stackrel{(i)}{\leq} \Theta \lp\tilde{\mu}\lp A\cap B_\infty(r)\rp\rp+ \Theta\lp\tilde{\mu}\lp B_\infty(r)^c\rp\rp\\
& \stackrel{(ii)}{\leq} 6C \overline{\Capa}_{\tilde{\mu}}(A\cap B_\infty(r))+\Theta\lp\tilde{\mu}\lp B_\infty(r)^c\rp\rp.\\
& \stackrel{(iii)}{\leq} 12C\lp\omega_T'\circ\omega_T^{-1}(r+1)\rp^2\Capa_{\tilde{\mu}}(A\cap B_\infty(r))+ 48dC\tilde{\mu}(B_\infty(r)^c)+\Theta\lp\tilde{\mu}\lp B_\infty(r)^c\rp\rp.\\
& \stackrel{(iv)}{\leq} 12C\lp\omega_T'\circ\omega_T^{-1}(r+1)\rp^2\Capa_{\tilde{\mu}}(A)+ \lp48dCT(1)+1\rp\Theta\lp\tilde{\mu}\lp B_\infty(r)^c\rp\rp,
\end{align*}
where (i) follows from the sub-additivity and the monotonicity of $\Theta$, (ii) from Lemma \ref{Cap-tilde}, (iii) from Lemma \ref{comp-Cap} and (iv) from the fact that the function $A\mapsto \Capa_{\tilde{\mu}}(A)$ is non decreasing and from the immediate inequality $x\leq T(1)\Theta(x)$ which holds for all $x\leq 1.$
Using Theorem \ref{concentration BLO}, it is not difficult to see that one can find $K\geq 1$ and $1\geq u_0>0$ such that $\tilde{\mu}\lp B_\infty(r)^c\rp\leq Ke^{-u_0r}$ for all $r\geq 0$. Thus $\Theta\lp\tilde{\mu}\lp B_\infty(r)^c\rp\rp\leq \Theta(Ke^{-u_0r})\leq K\Theta(e^{-u_0r})$ , where the last inequality follows from the sub-additivity of $\Theta$.
So, letting $a_1=12C\Capa_{\tilde{\mu}}(A)$, $a_2=(48dCT(1)+1)K$ and $t=\tilde{\mu}(A)$ and using the definitions of $\Theta$ and $\omega_T$, one has:
\[t\frac{1}{T\lp\frac{1}{\log(1+1/t)}\rp}\leq a_1 \frac{1}{T\lp\frac{1}{1+r}\rp}+a_2 e^{-u_0r}\frac{1}{T\lp\frac{1}{\log(1+e^{u_0r})}\rp}.\]
Using the inequality $1+u_0r\geq u_0(1+r)$ together with the sub-additivity property of the function $T$, one sees that
$\frac{1}{T\lp\frac{1}{1+r}\rp}\leq \frac{1}{u_0 T\lp\frac{1}{1+u_0r}\rp}.$
Thus the preceding inequality implies:
\[\forall v\geq 0,\quad t\frac{1}{T\lp\frac{1}{\log(1+1/t)}\rp}\leq \frac{a_1}{u_0} \frac{1}{T\lp\frac{1}{1+v}\rp}+a_2 e^{-v}\frac{1}{T\lp\frac{1}{\log(1+e^{v})}\rp}.\]
Now, observe that $(1+e^v)^3\geq 3e^v$, so $3\log(1+e^v)\geq \log(3)+v\geq 1+v$, and so $\frac{1}{T\lp\frac{1}{1+v}\rp}\leq\frac{1}{T\lp\frac{1}{3\log(1+e^{v})}\rp}\leq \frac{3}{T\lp\frac{1}{\log(1+e^{v})}\rp}$ where the last step follows from the sub-additivity property of $T$. So,
\[\forall v\geq 0,\quad t\frac{1}{T\lp\frac{1}{\log(1+1/t)}\rp}\leq \lp\frac{3a_1}{u_0} +a_2 e^{-v}\rp\frac{1}{T\lp\frac{1}{\log(1+e^{v})}\rp}.\]
Let $n\in \N^*$ ; taking $v=-n\log(t)$ in the preceding inequality gives:
\[t\frac{1}{T\lp\frac{1}{\log(1+1/t)}\rp}\leq \lp\frac{3a_1}{u_0} +a_2 t^n\rp\frac{1}{T\lp\frac{1}{\log(1+(1/t)^n)}\rp}.\]
Now, $(1+(1/t)^n)\leq (1+1/t)^n$, thus $\frac{1}{T\lp\frac{1}{\log(1+(1/t)^n)}\rp}\leq \frac{1}{T\lp\frac{1}{n\log(1+1/t)}\rp}\leq \frac{n}{T\lp\frac{1}{\log(1+1/t)}\rp}$ and consequently,
\[\frac{t}{n}\leq \lp\frac{3a_1}{u_0} +a_2 t^n\rp.\]
It is easy to check that if $n$ is sufficiently large, there is $m>0$ such that for all $t\in [0,1/2]$, one has $\frac{t}{n}-a_2 t^n\geq mt$. So $mt\leq \frac{3a_1}{u_0}$, that is to say \[\tilde{\mu}(A)\leq \frac{36C}{u_0m}\Capa_{\tilde{\mu}}(A).\]
A Capacity-measure inequality of this form is well known to imply the Poincar\'e inequality (see e.g \cite[Proposition 13 and Remark 20]{Barthe2006}).
\endproof

\bibliographystyle{plain}

\end{document}